\documentclass[12pt,a4paper]{article}
\usepackage{amssymb,times1,amsmath,amsfonts,amsthm,euscript}
\usepackage{hyperref}
\usepackage{graphicx}
\usepackage{mysec}
\usepackage{color}
\usepackage[normalem]{ulem}
\usepackage{cancel}
\newtheorem{theorem}{Theorem}[section]
\newtheorem{lemma}[theorem]{Lemma}

\theoremstyle{remark}
\newtheorem{example}{Example}[section]
\newtheorem{remark}{Remark}[section]

\theoremstyle{definition}

\newcommand{\PP}{\mathop{{\mathbb{P}}{}}\nolimits}
\newcommand{\EE}{\mathop{{\mathbb{E}}{}\myp}\nolimits}
\newcommand{\const}{\mathrm{const}}
\newcommand{\rd}{{\rm d}}
\newcommand{\re}{{\rm e}}
\newcommand{\ri}{{\rm i}}

\newcommand{\myp}{\mbox{$\:\!$}}
\newcommand{\mypp}{\mbox{$\;\!$}}
\newcommand{\myn}{\mbox{$\:\!\!$}}
\newcommand{\mynn}{\mbox{$\;\!\!$}}

\newcommand{\RR}{{\mathbb R}}

\newcommand{\NN}{{\mathbb N}}

\newcommand{\T}{T\kern-.06pc{}}

\newcommand{\calB}{{\mathcal B}}

\newcommand{\calF}{{\mathcal F}}

\makeatletter
\def\MR#1{\href{http://www.ams.org/mathscinet-getitem?mr=#1}{MR#1}}
\makeatother

\begin{document}

%\markboth{AUTHOR'S NAME (ALL CAPS)}{MANUSCRIPT TITLE (ALL CAPS)}

%% Place the running title of the paper with 40 letters or less in []
 %% and the full title of the paper in { }.
%\title[Archetypical equation: non-critical case] %Use the shortened version of the full title
%      {Bounded Continuous Solutions of the Archetypical Functional Equation: Analysis in the Noncritical Case}
\title{Analysis of the archetypal functional equation\\ in the non-critical
case\footnotetext{L.V.\,B.\ was partially supported by a Leverhulme
Research Fellowship. G.\myp{}D.\ received funding from the Israel
Science Foundation, Grant 35/10.
%Hospitality and
Support during various research visits provided by the Center for
Advanced Studies in Mathematics (Ben-Gurion University) and also by
the Center for Interdisciplinary Research (ZiF) and SFB\,701
(Bielefeld University) is acknowledged.}}

%%L.\,V.\,B.\ was partially supported by a Leverhulme Research
%%Fellowship, and
%% (2009--2011).

\author{Leonid V.~Bogachev\myp$^{\rm a}$, Gregory Derfel\myp$^{\rm b}$ and Stanislav A.~Molchanov\myp$^{\rm c}$}
\date{\small $^{\,\rm a}$\myp Department of Statistics, School of Mathematics, University of
Leeds,\\ Leeds LS2 9JT, UK. E-mail: {\tt
L.V.Bogachev@leeds.ac.uk}\\[.3pc]
$^{\,\rm b}$\myp Department of Mathematics, Ben Gurion University of
the Negev,\\
%P.O.\ Box 653,
Be'er Sheva 84105, Israel. E-mail: {\tt derfel@math.bgu.ac.il}\\[.3pc]
$^{\,\rm c}$\myp Department of Mathematics, University of North
Carolina at Charlotte,\\ Charlotte NC\,28223, USA. E-mail: {\tt
smolchan@uncc.edu}}

\maketitle
\begin{abstract}
We study the archetypal functional equation of the form
$y(x)=\iint_{\RR^2} y(a(x-b))\,\mu(\rd{}a,\rd{}b)$ ($x\in\RR$),
where $\mu$ is a probability measure on $\RR^2$; equivalently,
$y(x)=\EE\{y(\alpha\myp(x-\beta))\}$, where $\EE$ is expectation
with respect to the distribution $\mu$ of random coefficients
$(\alpha,\beta)$. Existence of non-trivial (i.e.\ non-constant)
bounded continuous solutions is governed by the value
$K:=\iint_{\RR^2}\ln\mynn|a|\,\mu(\rd{}a,\rd{}b)
=\EE\{\ln\mynn|\alpha|\}$; namely, under mild technical conditions
no such solutions exist whenever $K<0$, whereas if $K>0$ (and
$\alpha>0$) then there is a non-trivial solution constructed as the
distribution function of a certain random series representing a
self-similar measure associated with $(\alpha,\beta)$. Further
results are obtained in the supercritical case $K>0$, including
existence, uniqueness and a maximum principle. The case with
$\PP(\alpha<0)>0$ is drastically different from that with
$\alpha>0$; in particular, we prove that a bounded solution
$y(\cdot)$ possessing limits at $\pm\infty$ must be constant. The
proofs employ martingale techniques applied to the
martingale~$y(X_n)$, where $(X_n)$ is an associated Markov chain
with jumps of the form $x\rightsquigarrow\alpha\myp(x-\beta)$.


\medskip \noindent \emph{2010 MSC}:
%Mathematics Subject Classification}:
\,Primary: 39B05; \,Secondary: 34K06, 39A22, 60G42, 60J05.
\end{abstract}

% 39-XX Difference and functional equations
% 39Bxx Functional equations and inequalities [See also 30D05]
% 39B05 General
% 60J05 Discrete-time Markov processes on general state spaces
% 34K06 Linear functional-differential equations
% 39A22 Growth, boundedness, comparison of solutions
% 60G40 Stopping times; optimal stopping problems; gambling theory
% 60G42 Martingales with discrete parameter

\bigskip
\bigskip
\bigskip
\pagebreak
%\normalfont

%\author[L.~V.~Bogachev, G.~Derfel and S.~A.~Molchanov]{}

% Put your short thanks below. For long thanks/acknowlegements,
%please go to the last acknowlegments section.

%\bigskip

\section{Introduction}\label{sec0}
\subsection{The archetypal equation and main results}\label{sec:1.1} This paper
concerns the \emph{archetypal} functional equation with rescaled
argument \cite{BDM,Der} of the form
\begin{equation}\label{eq:eq0}
y(x)=\iint_{\RR^2} y(a(x-b))\,\mu(\rd{}a,\rd{}b),\qquad x\in\RR,
\end{equation}
where $\mu(\rd{}a,\rd{}b)$ is a probability measure on $\RR^2$. Due
to the normalization of the measure $\mu$ to unity, such an equation
is \textit{balanced} in that the total weighted contribution of the
(scaled) solution $y(\cdot)$ on the right-hand side of
\eqref{eq:eq0} is matched by the non-scaled input on the left-hand
side. The integral in \eqref{eq:eq0} has the meaning of expectation
with respect to a random vector $(\alpha,\beta)$ with distribution
$\PP\{(\alpha,\beta)\in \rd{a}\times\rd{b}\}=\mu(\rd{}a,\rd{}b)$;
thus, equation \eqref{eq:eq0} can be represented in the compact form
\begin{equation}\label{eq:arch}
y(x)=\EE\{y(\alpha\myp(x-\beta))\},\qquad x\in\RR.
\end{equation}

The equation \eqref{eq:eq0}--\eqref{eq:arch} is a rich source of
various equations specified by a suitable choice of the measure
$\mu$, which has motivated its name ``archetypal''~\cite{BDM}.
Examples include many well-known classes of equations with
rescaling, such as: equations in convolutions, e.g.\ the
\emph{Choquet--Deny equation} $y=y\star\sigma$ \cite{CD}; equations
for \emph{Hutchinson's self-similar measures} \cite{Hutchinson},
e.g.\ $y(x)=\tfrac12\mypp y(a\myp(x+1))+\tfrac12\mypp y(a\myp(x-1))$
($a>1$) arising in the \emph{Bernoulli convolutions problem}
\cite{Solomyak2}; \emph{two-scale} (\emph{refinement})
\emph{equations}\footnote{Compactly supported continuous solutions
of such equations play a crucial role in the theory of wavelets
\cite{D10,Strang}, and also in subdivision schemes and curve design
\cite{Cavaretta,DerDynLev}.} of the form $z(x)=a\sum_{i=1}^\ell
p_i\mypp z(a\myp(x-b_i))$ with $z(x):=y'(x)$ \cite{DL,DerDynLev},
exemplified by \emph{Schilling's equation}
$z(x)=a\myp\bigl(\frac{1}{4}\myp z(a\myp x - 1) +\frac{1}{2}\myp
z(a\myp x) + \frac{1}{4}\myp z(a\myp x + 1)\bigr)$ describing
spatially chaotic structures in amorphous materials
\cite{DerSchi,Schilling}; etc. Furthermore, as was observed by
Derfel \cite{Der}, the archetypal equation
\eqref{eq:eq0}--\eqref{eq:arch} also contains some important
functional-differential classes, including the (\emph{balanced})
\emph{pantograph equation}\footnote{Pantograph equation
$y'(x)=c_0\myp y(x)+c_1 y(\alpha x)$ dating back to Ockendon and
Tayler \cite{OT} arises in diverse areas, e.g.\ number theory,
astrophysics, radioactive decay, queues and risk theory, population
dynamics, medicine, quantum theory, stochastic games, etc.; for
general results and further bibliography on the pantograph equation,
see \cite{BDMO,BDM,BruntWake,DI,Iserles,KM}.} \cite{BDMO,BDM,Der}
\begin{equation}\label{eq:pant_gen}
y'(x)+y(x)=\sum\nolimits_{i} p_i\mypp y(a_i\myp x),\qquad
a_i,p_i>0,\quad \sum\nolimits_{i} p_i=1,
\end{equation}
and \emph{Rvachev's equation}\footnote{Its compactly supported
solutions are instrumental in approximation theory~\cite{Rvachev}.}
$z'(x)=2\bigl(z(2x+1)-z(2x-1)\bigr)$ \cite{Rvachev}. See an
extensive review of examples and applications of the archetypal
equation \eqref{eq:eq0}--\eqref{eq:arch} in Bogachev \emph{et al.}\
\cite{BDM}, together with further references therein.

Observing that any function $y(x)\equiv\const$ satisfies equations
\eqref{eq:eq0}--\eqref{eq:arch}, it is natural to investigate if
there are any \textit{non-trivial} (i.e.\ non-constant) bounded
continuous solutions. Such a question naturally arises in the
context of functional and functional-differential equations with
rescaling, where the possible existence of bounded solutions (e.g.\
periodic, almost periodic, compactly supported, etc.)\ is of major
interest in physical and other applications (see e.g.\
\cite{Cavaretta,Rvachev,Schilling,Strang}).

Investigation of the archetypal equation
\eqref{eq:eq0}--\eqref{eq:arch}, with a focus on bounded continuous
solutions (abbreviated below as \emph{b.c.-solutions}), was
initiated by Derfel \cite{Der} (in the case $\alpha>0$) who showed
that the problem crucially depends on the value
\begin{equation}\label{eq:K}
K:=\iint_{\RR^2}\ln\mynn|a|\,\mu(\rd{}a,\rd{}b)=\EE\{\ln\mynn|\alpha|\}.
\end{equation}
More precisely, if $K<0$ (\emph{subcritical case}) then, under some
mild technical conditions on the measure $\mu$, there are no
b.c.-solutions other than constants,\footnote{A similar result was
obtained earlier (via a different method) by Steinmetz and Volkmann
\cite{Steinmetz} for a special case of equation~\eqref{eq:arch},
\,$y(x) = p\mypp y(p\myp x-1) + q\mypp y(q\myp x + 1)$ \myp($p,q>0$,
$p+q=1$).} whereas if $K>0$ (\emph{supercritical case}) then a
non-trivial b.c.-solution does exist.

However, the \emph{critical case} $K=0$ was left open in~\cite{Der}.
Some recent progress was due to Bogachev \emph{et al.}\ \cite{BDMO}
who settled the problem for the balanced pantograph equation
\eqref{eq:pant_gen} by showing that if $K=\sum_{i} p_i\ln a_i=0$
then there are no non-trivial b.c.-solutions of~\eqref{eq:pant_gen}.
Recently (see \cite{BDM}) we proved the same result for a general
equation \eqref{eq:eq0}--\eqref{eq:arch} in the critical case
subject to an \emph{a priori} condition of uniform continuity
of~$y(\cdot)$, which is satisfied for a large class of examples
including~\eqref{eq:pant_gen}.

The focus of the present work is on the non-critical case $K\ne0$,
especially when $K>0$ with $\alpha$ possibly taking negative values,
aiming to obtain further results including existence, uniqueness and
a maximum principle. Under a slightly weaker moment condition on
$\beta$ as compared to \cite{Der} we establish the dichotomy of
non-existence \emph{vs.}\ existence of non-trivial b.c.-solutions in
the subcritical (${K<0}$) and supercritical ($K>0$) regimes,
respectively.

Let us stress though that in contrast to the subcritical case which
is insensitive to the sign of $\alpha$, for $K>0$ we are only able
to produce a non-trivial solution under the assumption that
$\alpha>0$ almost surely (a.s.). Such a solution is constructed,
with the help of results by Grintsevichyus \cite{Gr}, as the
distribution function $F_\Upsilon(x)=\PP(\Upsilon\le x)$ of the
random series $\Upsilon=\sum_{n=1}^\infty
\beta_n\prod_{i=1}^{n-1}\alpha_i^{-1}$ \strut{}representing a
self-similar measure associated with $(\alpha,\beta)$, where
$\{(\alpha_n,\beta_n)\}_{n\ge1}$ \strut{}are independent identically
distributed (i.i.d.)\ random pairs with distribution $\mu$ each.
This solution is unique (up to linear transformations) in the class
of functions with finite limits at $\pm\infty$
(Theorem~\ref{th:lim}\myp(a)), but the uniqueness in the class of
b.c.-solutions may fail to be true: we will present an example of
such a solution $y(\cdot)$ oscillating at $+\infty$ (see
Remark~\ref{rm:Kato}).

In the case $K>0$ with $\PP(\alpha<0)>0$, the function
$F_\Upsilon(\cdot)$ (which is still well defined) is no longer a
solution to the equation \eqref{eq:eq0}--\eqref{eq:arch}; e.g.\ if
$\alpha<0$ a.s.\ then $y=F_\Upsilon(x)$ satisfies another functional
equation, $y(x)=1-\EE\{y(\alpha\myp(x-\beta))\}$ (cf.\
\cite[Eq.\,(5)]{Gr}). Thus, the problem of existence remains largely
open here. More to the point, this case is completely different from
the purely positive case, $\alpha>0$ (a.s.); for instance, a
b.c.-solution $y(\cdot)$ with limits at $\pm\infty$ must be constant
(Theorem~\ref{th:lim}\myp(b)). This follows from Theorem
\ref{th:m+M<0} stating that the limits superior at $\pm\infty$
coincide (the same is true for the limits inferior). Heuristically,
this is a manifestation of ``mixing'' in \eqref{eq:arch} for (large)
positive and negative arguments of $y(\cdot)$ due to possible
negative values of $\alpha$. Note that Theorem \ref{th:m+M<0} is
proved with the help of the maximum principle of
Theorem~\ref{th:m+M>0}, which is of interest in its own right.

This analysis is complemented by uniqueness results in the class of
absolutely continuous (a.c.)\ solutions (using the Fourier transform
methods); here, boundedness is not assumed \emph{a priori}. Again,
we demonstrate a striking difference between the cases $\alpha>0$
(a.s.)\ and $\PP(\alpha<0)>0$ (see Theorems \ref{th:z1}
and~\ref{th:z2}, respectively).

Throughout the paper, it is assumed that
\begin{equation}\label{eq:as}
\text{(i)} \,\PP(\alpha\ne0)=1;\quad \text{(ii)}
\,\PP(|\alpha|\ne1)>0;\quad \text{(iii)} \ \forall\myp c\in\RR,\
\PP(\alpha\myp(c-\beta)=c)<1.
\end{equation}
Note that the remaining degenerate cases are treated in full detail
in~\cite{BDM}.

The rest of the paper is organized as follows. We start in
\S\ref{sec:Markov} by introducing an associated Markov chain $(X_n)$
with jumps of the form $x\rightsquigarrow\alpha\myp(x-\beta)$, and
also extend the iterated equation $y(x)=\EE_x\{y(X_n)\}$ to its
``optional stopping'' analog $y(x)=\EE_x\{y(X_\tau)\}$, where $\tau$
is a (random) stopping time and $\EE_x$ stands for the expectation
subject to the initial condition $X_0=x$. Suitable iterations of
such a kind will be instrumental. In \S\ref{sec:K<>0} we prove a
stronger version of the dichotomy between the cases $K<0$ and $K>0$
(the latter subject to $\alpha>0$). Finally, \S\ref{sec:super+}
contains further discussion of the supercritical case, as briefly
indicated above.

\section{Preliminaries}\label{sec:Markov}

\subsection{Associated Markov chain and harmonic functions}
The archetypal equation \eqref{eq:arch} admits an important
interpretation via an associated Markov chain $(X_n)$ on $\RR$
determined by the recursion
\begin{equation}\label{eq:Xn}
X_n=\alpha_n \myp(X_{n-1}-\beta_n) \quad (n\in\NN), \qquad
X_0=x\in\RR,
\end{equation}
where $\{(\alpha_n,\beta_n)\}_{n\ge1}$ are i.i.d.\ random pairs with
the same distribution as a generic copy $(\alpha,\beta)$. Transition
operator $\T$ of the Markov chain \eqref{eq:Xn} is given by
\begin{equation}\label{eq:generator}
\T f(x):= \EE_x\{f(X_1)\}\equiv\EE\{f(\alpha\myp(x-\beta))\},
\end{equation}
where the index $x$ indicates the initial condition $X_0=x$. A
function $f(\cdot)$ is called \emph{$T$-harmonic} (or simply
\emph{harmonic}) if $\T f=f$ (cf.\ \cite[p.\;40]{Revuz}); hence,
according to \eqref{eq:generator} solutions of equation
\eqref{eq:arch} are equivalently described as harmonic functions.

\subsection{Iterations and stopping times}
Equation \eqref{eq:arch} can be expressed as $y(x)=\EE_x\{y(X_1)\}$,
and by iteration $y(x)=\EE_x\{y(X_n)\}$ $(n\in\NN$). Explicitly,
\begin{gather}\label{eq:Y}
X_{n}=A_n x -D_n,\ \ \quad n\ge 0,\\
\label{eq:AY} A_{n}:=\prod_{i=1}^n \alpha_i\quad(A_0:=1),\ \ \quad
D_n:=\sum_{i=1}^{n}\beta_{i}\prod_{j=i}^n \alpha_{j}\quad(D_0:=0).
\end{gather}

For $n\in\NN_0:=\{0\}\cup\NN$, let $\calF_n:=\sigma\{X_i,\,i\le n\}$
be the $\sigma$-algebra generated by events $\{X_i\in B\}$ (with
Borel sets $B\in\calB(\RR)$); the increasing sequence
$(\calF_n)_{n\ge0}$ is referred to as the (natural) filtration of
$(X_n)$. A random variable $\tau$ with values in
${\NN\cup\{+\infty\}}$ is called a \emph{stopping time} with respect
to filtration $(\calF_n)$ if it is adapted to $(\calF_n)$ (i.e.\
$\{\tau=n\}\in\calF_n$, $n\in\NN_0$) and $\tau<\infty$ a.s. We shall
systematically use the following simple fact. (Note that the
continuity of $y(\cdot)$ is not required.)

\begin{lemma}\label{lm:stop1}
Let $\tau$ be a stopping time with respect to filtration
$\calF^\alpha_n:=\sigma\{\alpha_1,\dots,
\alpha_n\}\allowbreak\subset \calF_n$, $n\in\NN_0$. If\/ $y(\cdot)$
is a bounded $T$-harmonic function then
\begin{equation}\label{eq:stop_x}
y(x)=\EE_{x}\{y(X_{\tau})\},\qquad x\in\RR.
\end{equation}
\end{lemma}
\proof Clearly, $\tau$ is adapted to the filtration
$\mathcal{F}_n^{\alpha,\myp\beta}:=\sigma\{(\alpha_i,\beta_i),\,i\le
n\}\equiv \calF_n$. Using \eqref{eq:Xn} it is easy to check that
$\EE\{y(X_n)\mypp|\mypp\calF_{n-1}\}=y(X_{n-1})$ (a.s.),  and hence
$(y(X_n))$ is a \emph{martingale} \cite[p.\;43,
Proposition~1.8]{Revuz}. Since $y(\cdot)$ is bounded, formula
\eqref{eq:stop_x} now readily follows by Doob's Optional Stopping
Theorem \cite[pp.\ 485--486, Theorem~1 and Corollary]{Shiryaev}.
\endproof

\section{The subcritical ($K<0$) and supercritical
($K>0$) cases}\label{sec:K<>0} In the case $\alpha\ne0$ a.s.,
formula \eqref{eq:Y} can be rewritten in the form (cf.\
\eqref{eq:Y},~\eqref{eq:AY})
\begin{gather}\label{eq:Y*}
X_{n}=A_n (x -B_n),\ \ \quad n\ge 0,\\
\label{eq:AY*} A_{n}:=\prod_{i=1}^n \alpha_i\quad (A_0:=1),\ \ \quad
B_n:=D_nA^{-1}_n=\sum_{i=1}^{n}\beta_{i} A_{i-1}^{-1}\quad(B_0:=0).
\end{gather}

The following important result is due to Grintsevichyus
\cite[pp.\,164--165]{Gr}.
\begin{lemma}\label{lm:Grin}
Let assumption \eqref{eq:as} be in force, and also assume that
\begin{gather}\label{eq:12>0}
0<\EE\{\ln \mynn|\alpha|\}<\infty,
\qquad\EE\{\ln\max(|\beta|,1)\}<\infty.
\end{gather}
Then the random series
\begin{equation}\label{eq:Y-infty2}
\Upsilon:=\beta_1+\beta_2\mypp\alpha_1^{-1}+\beta_3\mypp\alpha_1^{-1}\myn\alpha_2^{-1}+\cdots
=\sum_{n=1}^\infty \beta_{n} A_{n-1}^{-1}
\end{equation}
converges a.s., and its distribution function
$F_\Upsilon(x):=\PP(\Upsilon\le x)$ is continuous on $\RR$.
\end{lemma}

\begin{remark}\label{rm:Grin}
The results in \cite{Gr} entail that $F_\Upsilon(\cdot)$ is either
a.c.\ or singularly continuous; a purely discrete case (with a
single atom!) arises if $\alpha\myp(c-\beta)=c$ (a.s.).
\end{remark}
Recall that the parameter $K$ is defined in~\eqref{eq:K}. The next
two results (for $K<0$ and $K>0$, respectively) were obtained by
Derfel \cite{Der} in the case $\alpha>0$ (a.s.)\ under a more
stringent condition $\EE\{|\beta|\}<\infty$; but his proofs
essentially remain valid in a more general situation as described
below.

\subsection{The subcritical case}

\begin{theorem}[$K<0$]\label{th:K<0}
Assume that the second integrability condition in \eqref{eq:12>0} is
fulfilled, but the first one is replaced by $-\infty<\EE\{\ln
\mynn|\alpha|\}<0$. Then any b.c.-solution of the archetypal
equation \eqref{eq:arch} is constant.
\end{theorem}
\proof Applying Lemma \ref{lm:stop1} with $\tau\equiv n\in\NN$, we
obtain (see \eqref{eq:Y},~\eqref{eq:AY})
\begin{equation}\label{eq:n-times}
y(x)=\EE\{y(A_n x-D_n)\},\qquad x\in\RR.
\end{equation}
Setting $D^\circ_n:=\sum_{i=1}^n \beta_i
A_i=\alpha_1(\beta_1+\beta_2\myp\alpha_2+\dots+\beta_n\alpha_2\cdots\alpha_n)$
(cf.~\eqref{eq:AY}), observe that the pair $(A_n,D_n)$ has the same
distribution as $(A_n,D^\circ_n)$, which is evident by reversing the
numbering $(\alpha_i,\beta_i)\mapsto(\alpha_{n-i+1},\beta_{n-i+1})$
\,($i=1,\dots,n$). Hence, equation \eqref{eq:n-times} can be
rewritten as
\begin{equation}\label{eq:n-times-tilde}
y(x)=\EE\{y(A_n x-D^\circ_n)\},\qquad x\in\RR.
\end{equation}
Due to Lemma \ref{lm:Grin} (with $\alpha_i^{-1}$ in place of
$\alpha_i$), $D^\circ_n$ converges a.s.\ as $n\to\infty$, say
$D^\circ_n\to\Upsilon^\circ$ (cf.~\eqref{eq:Y-infty2}). On the other
hand, $A_n\to 0$ a.s., since $\EE\{\ln\mynn|\alpha|\}<0$ and, by the
strong low of large numbers, $\ln\mynn |A_n|=\sum_{i=1}^n \ln\mynn
|\alpha_i|\to -\infty$ (a.s.). As a result, for each $x\in\RR$ we
have $A_nx-D^\circ_n\to-\Upsilon^\circ\!$ (a.s.). Since $y(\cdot)$
is bounded and continuous, one can apply Lebesgue's dominated
convergence theorem \cite[p.\,187, Theorem~3]{Shiryaev} and pass to
the limit in \eqref{eq:n-times-tilde}, yielding
$y(x)=\EE\{y(-\Upsilon^\circ)\}$; since the right-hand side does not
depend on $x$, it follows that $y(x)\equiv\const$.
\endproof

\subsection{Canonical solution in the supercritical case with $\alpha>0$}\label{sec:super}

The next theorem provides a non-trivial b.c.-solution to the
archetypal equation \eqref{eq:arch} in the case of \textit{positive}
$\alpha$. Recall that $\Upsilon$ is the random series
\eqref{eq:Y-infty2} and $F_\Upsilon(x)$ is its distribution function
(see Lemma~\ref{lm:Grin}).
\begin{theorem}[$K>0$]\label{th:K>0}
Suppose that assumption \eqref{eq:as} is satisfied, along with
conditions \eqref{eq:12>0}, and also assume that $\alpha>0$ a.s.
Then $y=F_\Upsilon(x)$ is a b.c.-solution of the archetypal
equation~\eqref{eq:arch}.
\end{theorem}

\proof Thanks to Lemma \ref{lm:Grin} we only have to verify that
$F_\Upsilon(x)$ satisfies \eqref{eq:arch}. Observe from
\eqref{eq:Y-infty2} that
$\Upsilon=\beta_1+\alpha_1^{-1}\myp\widetilde{\Upsilon}$, where
$\widetilde{\Upsilon}$ is independent of $(\alpha_1,\beta_1)$ and
has the same distribution as $\Upsilon$. Hence, we obtain (using
that $\alpha_1>0$ a.s.)
\begin{align*}
F_\Upsilon(x)&=\PP(\beta_1+\alpha_1^{-1}\myp\widetilde{\Upsilon} \le
x)= \PP(\widetilde{\Upsilon}\le \alpha_1(x-\beta_1))\\
&=\EE\bigl\{\PP\bigl(\widetilde{\Upsilon}\le
\alpha_1(x-\beta_1)\myp|\mypp
\alpha_1,\beta_1\bigr)\bigr\}=\EE\{F_\Upsilon(\alpha_1(x-\beta_1))\},
\end{align*}
that is, the function $y=F_\Upsilon(x)$ satisfies
equation~\eqref{eq:arch}.
\endproof

We will refer to $y=F_\Upsilon(x)$ as the \textit{canonical
solution} of equation~\eqref{eq:arch}.

\begin{remark}
For some concrete equations with $\alpha\equiv\const>1$,
b.c.-solutions different from the canonical one may be constructed
(see Remark~\ref{rm:Kato}).
\end{remark}

\begin{remark}
To the best of our knowledge, no non-trivial b.c.-solutions of
equation \eqref{eq:arch} are known if $\PP(\alpha<0)>0$ except in
the special case $|\alpha|\equiv 1$ (see \cite[Theorem
2.2\myp(b-ii)]{BDM}).
\end{remark}

\section{Further results in the supercritical case}\label{sec:super+}
\subsection{Bounds coming from infinity}
The next result is akin to the maximum principle for the usual
harmonic functions. The continuity of $y(\cdot)$ is not presumed.

\begin{theorem}[Maximum Principle]\label{th:m+M>0}
Suppose that assumption \eqref{eq:as} is satisfied, along with
conditions \eqref{eq:12>0}. Let\/ $y(\cdot)$ be a bounded solution
of \eqref{eq:arch}, and denote
\begin{equation}\label{eq:lim}
m^\pm:=\liminf_{x\to\pm\infty} y(x),\qquad
M^\pm:=\limsup_{x\to\pm\infty} y(x),
\end{equation}
where the same $+$ or $-$ sign should be chosen on both sides of\/
each equality. Then
\begin{equation}\label{eq:m+M}
m\le  y(x)\le M,\qquad x\in\RR,
\end{equation}
where $m:=\min\{m^{+}\myn,m^{-}\}$, \,$M:=\max\{M^{+}\myn,M^{-}\}$.
\end{theorem}

\proof Applying Lemma \ref{lm:stop1} with $\tau\equiv n\in\NN$, for
any $x\in\RR$ we obtain
\begin{equation}\label{eq:iter-n}
y(x)=\EE\{y(A_n(x-B_n))\},
\end{equation}
where $A_n=\prod_{i=1}^n \alpha_i$ and $B_n=\sum_{i=1}^n \beta_{i}
A_{i-1}^{-1}$ (see \eqref{eq:Y*},~\eqref{eq:AY*}). By Lemma
\ref{lm:Grin}, the limiting random variable
$\Upsilon=\lim_{n\to\infty} B_n$ is continuous, hence
$\lim_{n\to\infty} (x-B_n)=x-\Upsilon\ne0$ (a.s.). Combined with
$|A_n|\to\infty$ a.s.\ (which follows by the strong law of large
numbers due to the first moment condition in \eqref{eq:12>0}, cf.\
the proof of Theorem~\ref{th:K<0}), this implies that
$|A_n(x-B_n)|\to\infty$ (a.s.). Hence, Fatou's lemma \cite[p.\,187,
Theorem~2]{Shiryaev} applied to equation \eqref{eq:iter-n} yields
\begin{equation*}
y(x) \le \EE\!\left\{\limsup_{n\to\infty}y(A_n(x-B_n))\right\} \le
\max\{M^{+}\myn,M^{-}\}=M,
\end{equation*}
which proves the upper bound in \eqref{eq:m+M}. The lower bound
follows similarly.
\endproof

The case where $\alpha$ may take on negative values has an
interesting general property as follows (note that conditions
\eqref{eq:12>0} are not needed here).
\begin{theorem}\label{th:m+M<0}
Suppose that $q:=\PP(\alpha<0)>0$, and let $y(x)$ be a bounded
solution of~\eqref{eq:arch}. Then, in the notation \eqref{eq:lim},
we have
\begin{equation}\label{eq:inf+sup}
m^{-}\!=m^{+},\qquad M^{-}\!=M^{+}.
\end{equation}
\end{theorem}
\proof By Fatou's lemma applied to equation \eqref{eq:arch} we get
\begin{align}
M^{+}=\limsup_{x\to+\infty}\,y(x) &\le
\EE\!\left\{\limsup_{x\to+\infty} \,y(\alpha\myp(x-\beta))\right\}
\label{eq:M<M} \le M^{+} (1-q)+M^{-} q.
\end{align}
Since $q>0$, \eqref{eq:M<M} implies that $M^{+}\!\le M^{-}$. By
symmetry, the opposite inequality is also true, hence $M^{-}\!=
M^{+}$. The first equality in \eqref{eq:inf+sup} is proved
similarly.
\endproof

\subsection{Uniqueness for solutions with limits at infinity}
We can now prove the following uniqueness result (again, the
continuity of solutions is not presumed). Note that the cases
$\alpha>0$ (a.s.) and $\PP(\alpha<0)>0$ are drastically different.

\begin{theorem}\label{th:lim}
Let assumption \eqref{eq:as} be in force, along with
conditions~\eqref{eq:12>0}. Let $y(\cdot)$ be a bounded solution of
\eqref{eq:arch} such that the limits $L^\pm:=\lim_{x\to\pm\infty}
y(x)$ exist.
\begin{itemize}
\item[\rm (a)] Suppose that\/ $\PP(\alpha>0)=1$. Then
$y(\cdot)$ coincides, up to an affine transformation, with the
canonical solution $F_{\Upsilon}(\cdot)$ (see Theorem~\ref{th:K>0});
specifically,
\begin{equation}\label{eq:L0L1}
y(x)=(L^{+}-L^{-})\mypp F_\Upsilon(x) + L^{-}, \qquad x\in\RR.
\end{equation}
In particular, $y(\cdot)$ must be everywhere continuous.

\smallskip
\item[\rm (b)] If\/ $\PP(\alpha<0)>0$ then $y(x)\equiv
\const$.
\end{itemize}
\end{theorem}

\proof (a) Denote the right-hand side of \eqref{eq:L0L1} by
$y_*(x)$. By linearity of \eqref{eq:arch} and according to Theorem
\ref{th:K>0}, $y_*(x)$ satisfies equation \eqref{eq:arch}, and it
has the same limits $L^\pm$ at $\pm\infty$ as the solution $y(x)$.
Hence, $y(x)-y_*(x)$ is also a solution, with zero limits at
$\pm\infty$. But Theorem \ref{th:m+M>0} then implies that
$y(x)-y_*(x)\equiv 0$.

\smallskip
(b) Theorem \ref{th:m+M<0} implies that $L^-\myn=L^+=:L$, hence by
the bound \eqref{eq:m+M} of Theorem \ref{th:m+M>0} we have $L\le
y(x)\le L$, i.e.\ $y(x)\equiv L=\const$.
\endproof

\begin{remark}
In the case  $\PP(\alpha<0)>0$, Theorem \ref{th:lim}\myp(b) holds
true if just one of the limits $L^{\pm}$ is assumed (due to
\eqref{eq:inf+sup}, the other limit exists automatically).
\end{remark}

\begin{remark}\label{rm:Kato}
Kato and McLeod \cite[p.\;923, Theorem~9\myp(iii)]{KM} showed
\emph{inter alia} that the pantograph equation $y'(x)+y(x)=y(\alpha
x)$ with $\alpha=\const>1$ has a family of $C^\infty$-solutions on
the half-line $x\in[0,\infty)$ such that $y(x)={\rm g}(\ln
x/\ln\alpha)+O(x^{-\theta})$ as $x\to+\infty$, where ${\rm
g}(\cdot)$ is any $1$-periodic function, H\"older continuous with
exponent $0<\theta\le1$. Noting from the equation that $y'(0)=0$,
such solutions can be extended to the entire line $\RR$ by defining
$y(x):=y(0)$ for all $x<0$. It is known (see \cite{BDM,Der}) that
$y(\cdot)$ automatically satisfies the archetypal equation
\eqref{eq:arch} (with the same $\alpha>1$ and exponentially
distributed $\beta$), thus furnishing an example of (a family of)
bounded continuous (even smooth) solutions that do not have limit at
$+\infty$.
\end{remark}

\subsection{Uniqueness via Fourier transform}
Here, we obtain uniqueness results in the class of a.c.\ solutions
with integrable derivative. In what follows, abbreviation ``a.e.''
stands for ``almost everywhere'' (with respect to Lebesgue measure
on $\RR$). Note that boundedness of solutions is not presumed. It is
convenient to state and prove these results separately for positive
and negative~$\alpha$ (see Theorems \ref{th:z1} and \ref{th:z2},
respectively). Recall that $\Upsilon$ is the random
series~\eqref{eq:Y-infty2}.

\begin{theorem}\label{th:z1}
Let assumption \eqref{eq:as} be satisfied, together with
conditions~\eqref{eq:12>0}.
\begin{itemize}
\item[\rm (a)] Let\/ $\alpha>0$ a.s., and assume that a solution
$y(\cdot)$ of equation \eqref{eq:arch} is a.e.\ differentiable, with
$z(x):=y'(x)\in L^1(\RR)$. Then $z(\cdot)$ is determined uniquely
\textup{(}a.e.\textup{)}\ up to a multiplicative factor, with
Fourier transform given by
\begin{equation}\label{eq:Fourier_infty}
\hat{z}(s)=c_1\EE\{\re^{\myp\ri s\Upsilon}\}\ \quad (s\in\RR),\qquad
c_1:=\hat{z}(0)\in\RR.
\end{equation}
\item[\rm (b)] If
$y(\cdot)$ is also a.c.\ then it coincides, up to an affine
transformation, with the canonical solution $F_{\Upsilon}(\cdot)$
(see Theorem~\ref{th:K>0}), i.e.\ there are $c_0,c_1\in\RR$ such
that
\begin{equation}\label{eq:c0c1}
y(x)=c_0+c_1 F_\Upsilon(x),\qquad x\in\RR.
\end{equation}
\end{itemize}
\end{theorem}
\proof (a) Differentiation of \eqref{eq:arch} shows that
$z(x):=y'(x)$ satisfies a.e.\ the equation
\begin{equation}\label{eq:z}
z(x)=\EE\{\alpha\myp z(\alpha\myp(x-\beta))\}.
\end{equation}
Let $\hat{z}(s):=\int_\RR \re^{\myp\ri sx}z(x)\,\rd{x}$ be the
Fourier transform of the function $z\in L^1(\RR)$, hence
$\hat{z}(\cdot)$ is bounded and continuous on $\RR$, with the
sup-norm $\|\hat{z}\|\le \int_\RR |z(x)|\,\rd{x}<\infty$.
Multiplying \eqref{eq:z} by $\re^{\myp\ri sx}$ and integrating over
$x\in\RR$, we see, using Fubini's theorem and the substitution
$t=\alpha\myp(x-\beta)$, that $\hat{z}(\cdot)$ satisfies the
equation
\begin{equation}\label{eq:Fourier}
\hat{z}(s)=\EE\{\re^{\myp\ri s\beta}\myp
\hat{z}(\alpha^{-1}s)\},\qquad s\in\RR.
\end{equation}
Iterating \eqref{eq:Fourier} $n\ge1$ times we get (see the
notation~\eqref{eq:AY*})
\begin{equation}\label{eq:Fourier_n}
\hat{z}(s)=\EE\bigl\{\re^{\myp\ri
sB_n}\mypp\hat{z}(A_n^{-1}s)\bigr\},\qquad s\in\RR.
\end{equation}
Note that $\EE\{\ln\mynn|\alpha^{-1}|\} \in (-\infty,0)$, hence
$A_n^{-1}\to0$ a.s.\ (see the proof of Theorem \ref{th:K<0});
besides, $B_n\to\Upsilon$ a.s.\ by Lemma \ref{lm:Grin}. Thus,
passing to the limit in \eqref{eq:Fourier_n} (by dominated
convergence) and recalling that $\hat{z}(\cdot)$ is continuous, we
obtain~\eqref{eq:Fourier_infty}.

(b) To identify $z(\cdot)$ from its Fourier transform
\eqref{eq:Fourier_infty}, it is convenient to integrate both parts
of equation \eqref{eq:Fourier_infty} against a suitable class of
test functions. Consider the Schwartz space $\mathcal{S}(\RR)$ of
smooth functions $\varphi(x)$ with finite support and such that
their Fourier transform $\widehat{\varphi}(s)=\int_\RR \re^{\myp \ri
sx}\myp\varphi(x)\,\rd x$ is integrable; by the inversion formula,
$\varphi(x)=(2\pi)^{-1}\!\int_\RR \re^{-\myp\ri
sx}\:\widehat{\varphi}(s)\,\rd s$. With this at hand, we can write
\begin{align}
\int_{\RR} \hat{z}(s)\,\widehat{\varphi}(s)\,\rd{s} \label{Pars1}
&=\int_{\RR} \left(\int_\RR \re^{\myp\ri
sx}\:\widehat{\varphi}(s)\,\rd{s}\right)z(x)\,\rd{x}=2\pi\int_{\RR}
\varphi(-x)\, z(x)\,\rd x.
\end{align}
Similarly,
\begin{align}
\notag \int_{\RR} \EE\{\re^{\myp\ri
s\Upsilon}\}\,\widehat{\varphi}(s)\,\rd{s}&=\int_{\RR}
\left(\int_\RR \re^{\myp\ri
sx} \,\rd F_\Upsilon(x)\right)\widehat{\varphi}(s)\,\rd{s}\\
\label{Pars2} &=\int_{\RR} \left(\int_\RR \re^{\myp\ri
sx}\:\widehat{\varphi}(s)\,\rd{s}\right)\rd F(x)=2\pi \int_{\RR}
\varphi(-x)\,\rd F_\Upsilon(x).
\end{align}
Thus, thanks to equation \eqref{eq:Fourier_infty}, from
\eqref{Pars1} and \eqref{Pars2} we obtain
\begin{equation}\label{eq:1=2}
\int_{\RR} \varphi(-x)\, z(x)\,\rd x= c_1\int_{\RR} \varphi(-x)\,\rd
F_\Upsilon(x),\qquad \varphi\in\mathcal{S}(\RR).
\end{equation}
Since $\mathcal{S}(\RR)$ is dense in both $L^1(\RR;\,z(x)\,\rd x)$
and $L^1(\RR;\rd F_\Upsilon(x))$, equation \eqref{eq:1=2} extends to
indicator functions of any intervals, yielding (by the continuity of
$F_\Upsilon(\cdot)$)
$$
y(x)-y(0)=\int_{0}^{x} z(u)\,\rd u=c_1
\{F_\Upsilon(x)-F_\Upsilon(0)\},\qquad x\in\RR,
$$
which is reduced to \eqref{eq:c0c1} by setting $c_0:=y(0)- c_1
F_\Upsilon(0)$.
\endproof

\begin{remark}
The result of Theorem \ref{th:z1} was obtained by Daubechies and
Lagarias \cite[p.\,1392, Theorem~2.1(b)]{DL} in a particular case
with $\alpha\equiv\const>1$ and discrete~$\beta$.
\end{remark}

\begin{remark}
Uniqueness (up to a multiplicative factor) of b.c.-solutions of
equation \eqref{eq:Fourier} was proved by Grintsevichyus
\cite[p.\,165, Proposition~l]{Gr}.
\end{remark}

%\begin{remark}
%Under the hypotheses of Theorem~\ref{th:z1}\myp(b), the function
%$F_\Upsilon(\cdot)$ must be a.c.\ (cf.\ Remark~\ref{rm:Grin}). In
%other words, if $F_\Upsilon(\cdot)$ happens to be singular then
%there are no a.c.\ solutions of \eqref{eq:arch} with integrable
%derivative;
%% see \cite[Theorem~2,p.\,165]{Gr}
%e.g., this is the case for $\alpha\equiv\const>2$ and $\beta=\pm1$
%with probabilities $\frac12$ (see~\cite{Kershner}).
%\end{remark}

\begin{example}
De Rham's function (see \cite[pp.\ 1403--1405]{DL} is a continuous
(but nowhere differentiable) even solution of the difference
equation
$$
\phi(x)=\phi(3x)+\tfrac13\bigl(\phi(3x + 1) +\phi(3x-1)\bigr)
+\tfrac23\bigl(\phi(3x + 2) +\phi(3x - 2)\bigr).
$$
Then $y(x):=\int_0^x \phi(u)\,\rd{u}$ is an odd function of class
$C^1(\RR)$ satisfying
$$
y(x)=\tfrac13\mypp y(3x)+\tfrac19\bigl(y(3x + 1) +y(3x-1)\bigr)
+\tfrac29\bigl(y(3x + 2) +y(3x - 2)\bigr),
$$
which is an archetypal equation with $\alpha\equiv 3$ and $\beta$
taking values $0,-\tfrac13,\tfrac13,-\tfrac23,\tfrac23$ with
probabilities $\tfrac13,\tfrac19,\tfrac19,\tfrac29,\tfrac29$,
respectively. Now, according to Theorem \ref{th:z1} the solution
$y(\cdot)$ is an affine version of the distribution function
$F_\Upsilon(\cdot)$, the latter thus being automatically a.c.\ and,
moreover, in $C^1(\RR)$; in turn, it follows that de Rham's function
$\phi(\cdot)$ is proportional to the probability density of
$\Upsilon$ (see~\eqref{eq:Y-infty2}).
\end{example}

A counterpart of Theorem \ref{th:z1} for $\alpha$ with possible
negative values is strikingly different (cf.\ Theorem~\ref{th:lim}).
\begin{theorem}\label{th:z2}
Let $q:=\PP(\alpha<0)>0$, and let a solution $y(\cdot)$ be a.e.\
differentiable, with $y'\in L^1(\RR)$. Then $y'=0$ a.e. If in
addition $y(\cdot)$ is a.c.\ then $y\equiv \const$.
\end{theorem}
\proof The random time $\tau_-:=\inf\{n\ge 1\colon A_n<0\}$ is
adapted to the filtration $\calF^\alpha_n$ and has geometric
distribution, $\PP(\tau_{-}\!=n)=\allowbreak(1-q)^{n-1}q$
\,($n\ge1$). Hence, $\tau_{-}\!<\infty$ a.s.\ and
$\EE\{\tau_{-}\}=q^{-1}<\infty$. Applying Lemma \ref{lm:stop1}, we
obtain the equation
\begin{equation}\label{eq:tilde-eq}
y(x)=\EE\{y(\tilde\alpha\myp(x-\tilde\beta)\},\qquad x\in\RR,
\end{equation}
where $\tilde\alpha:=A_{\tau_{-}}<0$, \,$\tilde\beta:=B_{\tau_{-}}$
(cf.\ \eqref{eq:Y*},~\eqref{eq:AY*}).

Let us first verify that $\tilde\alpha$, $\tilde\beta$ satisfy the
moment conditions \eqref{eq:12>0}. Indeed, noting that
$\ln\mynn|\tilde\alpha| =\sum_{i=1}^{\tau_{-}}\ln\mynn|\alpha_i|$
and $\EE\{\tau_{-}\}=q^{-1}<\infty$, by Wald's identity
\cite[p.\;488, Theorem~3]{Shiryaev} we obtain, using the first
condition in~\eqref{eq:12>0},
\begin{equation}\label{eq:*} \EE\{\ln\mynn|\tilde\alpha|\}=\EE\{\tau_{-}\}\cdot
\EE\{\ln\mynn|\alpha|\}\in(0,\infty).
\end{equation}
Recalling \eqref{eq:AY*} and denoting $a\vee b:=\max\{a,b\}$,
\,$a\wedge b:=\min\{a,b\}$, we have
\begin{align*}
|\tilde\beta|\le \sum_{i=1}^{\tau_{-}} \frac{|\beta_i|}{|A_{i-1}|}
&\le\prod_{i=1}^{\tau_{-}} (|\beta_i|\vee 1)\cdot
\sum_{i=1}^{\tau_{-}}\frac{1}{|A_{i-1}|} \prod_{i=1}^{\tau_{-}}
(|\beta_i|\vee 1)\cdot
\tau_{-}\prod_{i=1}^{\tau_{-}}\frac{1}{|\alpha_i|\wedge 1}.
\end{align*}
The right-hand side is not less than $1$, hence the same bound holds
for $|\tilde\beta|\vee 1$ and
\begin{equation}\label{eq:prod-sum2}
\ln\mynn(|\tilde\beta|\vee 1)\le \sum_{i=1}^{\tau_{-}}
\ln\mynn(|\beta_i|\vee 1)+\ln
\mynn(\tau_{-})-\sum_{i=1}^{\tau_{-}}\ln\mynn (|\alpha_{i}|\wedge
1).
\end{equation}
Again applying Wald's identity and using conditions \eqref{eq:12>0},
we get from~\eqref{eq:prod-sum2}
\begin{equation*}
\EE\{\ln\mynn(|\tilde\beta|\vee 1)\}\le
\EE\{\tau_{-}\}\cdot\Bigl(\EE\{\ln\mynn(|\beta|\vee
1)\}+1-\EE\{\ln\mynn(|\alpha|\wedge 1)\}\Bigr)<\infty.
\end{equation*}

Now we can apply to \eqref{eq:tilde-eq} the method used in the proof
of Theorem \ref{th:z1}. More specifically, a differentiated version
of \eqref{eq:tilde-eq}, for $z(x):=y'(x)$, reads (cf.~\eqref{eq:z})
\begin{equation*}
z(x)=\EE\{\tilde{\alpha}\mypp z(\tilde{\alpha}\myp
(x-\tilde{\beta}))\}\qquad \text{(a.e.)}.
\end{equation*}
However, here $\tilde\alpha<0$ (a.s.), so the Fourier transform
$\hat{z}(s)$ now satisfies (cf.~\eqref{eq:Fourier})
\begin{equation*}
\hat{z}(s)=-\EE\{\re^{\myp\ri s\tilde{\beta}}\myp
\hat{z}(\tilde{\alpha}^{-1}s)\},\qquad s\in\RR.
\end{equation*}
Iterating as before, we obtain for each $n\in \NN$
\begin{equation}\label{eq:Fourier_n_tilde}
\hat{z}(s)=(-1)^n\EE\{\re^{\myp\ri
s\widetilde{\Upsilon}_n}\mypp\hat{z}(\tilde{A}_n^{-1}s)\},\qquad
s\in\RR,
\end{equation}
where due to \eqref{eq:*} we have a.s. $\tilde{A}_n^{-1}
=\prod_{i=1}^n \tilde{\alpha}_i^{-1} \to0$,
\,$\widetilde{\Upsilon}_n=\sum_{i=1}^n \tilde{\beta}_{i}
\tilde{A}_{i-1}^{-1}\to \widetilde{\Upsilon}$. Hence, the
expectation in \eqref{eq:Fourier_n_tilde} converges to
$\hat{z}(0)\EE\{\re^{\myp\ri s\widetilde{\Upsilon}}\}$; however, due
to the sign alternation the limit of \eqref{eq:Fourier_n_tilde} does
not exist unless $\hat{z}(0)=0$, in which case $\hat{z}(s)=0$ for
all $s\in\RR$. By the uniqueness theorem for the Fourier transform,
this implies that $z(x)= y'(x)\equiv 0$ a.e. Finally, if $y(\cdot)$
is a.c.\ then it follows that $y(x)\equiv\const$.
\endproof

\begin{remark}
The last statement (i.e.\ under the a.c.-condition) of each of
Theorems \ref{th:z1} and \ref{th:z2} can be easily deduced by
Theorem~\ref{th:lim}. Indeed, since the derivative $y'(\cdot)$ is
a.c.\ and in $L^1(\RR)$, by the Newton--Leibniz formula we have
$$
y(x)=y(0)+\int_0^x y'(u)\,\rd{u}\to  y(0)+\int_0^{\pm\infty}
y'(u)\,\rd{u}\qquad (x\to\pm\infty).
$$
Thus, the limits of $y(x)$ at $\pm\infty$ exist, and the rest
immediately follows from Theorem~\ref{th:lim}. However, the
uniqueness results for the derivative $y'$, contained in Theorems
\ref{th:z1} and \ref{th:z2}, cannot be obtained along these lines.
\end{remark}

\subsubsection*{Acknowledgments.}
The authors are grateful to John Ockendon and Anatoly Vershik for
stimulating discussions.

\end{document}